\documentclass[12pt, reqno]{amsart}
\usepackage{latexsym}
\usepackage{mathtools}
\usepackage{amsmath}
\usepackage{amsthm}
\usepackage{mathrsfs}
\usepackage{lmodern}
\usepackage[T1]{fontenc}
\usepackage{textcomp}
\usepackage[utf8]{inputenc}
\usepackage{color}
\usepackage{amssymb}
\usepackage{enumerate}
\usepackage[abbrev]{amsrefs}
\usepackage{indentfirst}
\usepackage{graphicx}
\usepackage{bmpsize}
\usepackage{hyperref}
\usepackage{enumitem}
\usepackage{enumerate}
\usepackage{url}
\usepackage{autobreak}
\usepackage[justification=justified]{caption}
\usepackage{bm}
\usepackage{here}
\usepackage[labelformat=empty,subrefformat=parens]{subcaption}
\renewcommand{\thefootnote}{} %footnote counter

\theoremstyle{plain} %text of this environment is typesetted in italics
\newtheorem{theorem}{\indent\sc Theorem}[section]
\newtheorem{lemma}[theorem]{\indent\sc Lemma}
\newtheorem{corollary}[theorem]{\indent\sc Corollary}
\newtheorem{proposition}[theorem]{\indent\sc Proposition}

\newtheorem{subpath}{\indent\sc Subpath}
\theoremstyle{definition} %text of this environment is typesetted in roman letters
\newtheorem{definition}[theorem]{\indent\sc Definition}
\newtheorem{remark}[theorem]{\indent\sc Remark}

\newenvironment{sproof}{%
  \proof}{\endproof}
\makeatletter
  
  \@addtoreset{equation}{section}
\makeatother
\newcommand{\Ca}{\mathfrak{C}}

\allowdisplaybreaks[3] %改ページを数式に許す

\newcommand\cinput[2]{\lower#1pt\hbox{\input{#2}}}

\begin{document}

\title{Divergence functions of higher-dimensional Thompson's groups}
\author{Yuya Kodama}
\date{}
\renewcommand{\thefootnote}{\arabic{footnote}}  %number
\setcounter{footnote}{0} %footnote counter
%\thanks{The author was supported by JST, the establishment of university fellowships towards the creation of science technology innovation, Grant Number JPMJFS2139.}

\begin{abstract}
We prove that higher-dimensional Thompson's groups have linear divergence functions. 
By the work of Dru{\c{t}}u, Mozes, and Sapir, this implies none of the asymptotic cones of $nV$ has a cut-point. 
\end{abstract}

\maketitle

\section{Introduction}
%%%%%%%%%%%%%%%%%%%%%%%%%%%%%%%%%%%%%%%%%
Thompson's groups $F$, $T$, and $V$ are finitely presented infinite groups defined by Richard Thompson in the 1960s. 
They are all known to be mysterious groups. 
For example, $T$ and $V$ are the first examples of finitely presented,  infinite, and simple, and it is known that the amenability of $F$ is a difficult open problem. 
Because they have several unpredictable properties, by focusing on such properties, many ``generalized'' Thompson's groups were also defined. 
Higher-dimensional Thompson's groups, denoted by $2V, 3V, \cdots$, are also such groups defined by Brin \cite{MR2112673}. 
The group $V$ acts on the Cantor set $\Ca$, and the group $nV$ acts on the powers of the Cantor set $\Ca^n$. 
It is known that $nV$ is also finitely presented \cite{MR2114568, hennig2012presentations} and simple \cite{MR2112673, brin2010baker}. 
In addition, it was shown that for $n, m \in \mathbb{Z}_{>0}$, the group $nV$ is isomorphic to $mV$ if and only if $n=m$ holds \cite{bleak2010family}. 

In 2018, Golan and Sapir showed that the divergence functions of $F$, $T$, and $V$ are linear \cite{golan2019divergence}. 
This function was first mentioned by Gromov\cite{zbMATH00437296}, and later, Gersten gave the formal definition \cite{gersten1994quadratic} as a quasi-isometric invariant of geodesic metric spaces. 
Roughly speaking, the order of the function indicates whether the Cayley graphs of the group are ``close'' to the Euclidean or hyperbolic spaces. 
In fact, the orders of the functions of the direct sums of the infinite cyclic groups $\mathbb{Z}^2, \mathbb{Z}^3, \cdots$ are linear, and it is known that the orders of the functions of hyperbolic groups are at least exponential \cite{bridson2013metric}. 
In \cite{golan2019divergence}, they asked whether their proof could be extended to generalized Thompson's groups. 
In recent years, similar results have been obtained for some groups by extending the original arguments \cite{kodama2023divergence, sheng2024divergence, lucy2024divergence}. 
For the recent results of the functions of groups other than generalized Thompson's groups, see \cite{issini2023linear}. 

In this paper, we also extend the argument given by Golan and Sapir. 
That is, we show the following:  
\begin{theorem}\label{theorem_divergence_nV}
Higher-dimensional Thompson's groups have linear divergence functions. 
\end{theorem}

This paper is organized as follows: 
in Section \ref{section_preliminary}, we summarize the definition of higher-dimensional Thompson's group $2V$ and the notion of divergence functions. 
In Section \ref{section_main}, we first prepare estimates of the word length of elements in $2V$ to ensure that the condition required by the definition of the divergence function is satisfied. 
Subsequently, we construct a ``good'' path from any $g \in 2V$ to a specific element in $2V$, where the element is determined only from the word length of $g$. 
It should be noted that we assume $n=2$ in most of this paper; however, all the proofs can be generalized to $n$ with the appropriate modifications. 
%%%%%%%%%%%%%%%%%%%%%%%%%%%%%%%%%%%%%%%%%
\section{Preliminaries} \label{section_preliminary}
%\subsection{Thompson's group $V$}
%%%%%%%%%%%%%%%%%%%%%%%%%%%%%%%%%%%%%%%%%
\subsection{Definition of higher-dimensional Thompson's groups}
In this section, we define only higher dimensional Thompson's group $2V$. 
For the readers unfamiliar with Thompson's groups, see \cite{cannon1996introductory}. 
\subsubsection{Homeomorphisms on the direct product of Cantor sets}
Let $\Ca$ be the Cantor set $\{0, 1\}\times \{0, 1\} \times \cdots$. 
For a finite word $w$ on $\{0, 1\}$ and a finite or infinite word $\zeta$ on $\{0, 1\}$, let $w\zeta$ denote their concatenation. 
Following \cite{MR2112673}, we will describe partitions of $\Ca^2$ by using subdivisions of the unit square $[0, 1]^2$. 
Subsequently, we will define homeomorphisms from $\Ca$ to itself that are obtained from such partitions. 
%We first recall the definition of Thompson's group $V$ and then review the definition of higher-dimensional Thompson's group as a generalization of $V$. 

\begin{figure}[tbp]
\begin{center}
\includegraphics[width=0.7\linewidth]{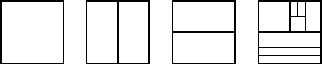}
\end{center}
\caption{The trivial pattern, the pattern once divided vertically, the pattern once divided horizontally, and a pattern. }
\label{Fig_patterns}
\end{figure}

Firstly, we call $[0, 1]^2$ itself \textit{trivial pattern}. 
Let us consider a rectangle $[a_1, a_2]\times [b_1, b_2] \subset [0, 1]^2$. 
By dividing this rectangle in half, we obtain two new rectangles. 
The way of obtaining rectangles $[a_1, (a_1+a_2)/2] \times [b_1, b_2]$ and $[(a_1+a_2)/2, a_2] \times [b_1, b_2]$ is called \textit{vertical subdivision}, and the way of obtaining rectangles $[a_1, a_2]\times [b_1, (b_1+b_2)/2]$ and $[a_1, a_2] \times [(b_1+b_2)/2, b_2]$ is called \textit{horizontal subdivision}. 
A \textit{pattern} is defined as a finite set of rectangles in $[0, 1]^2$ obtained from the trivial pattern by applying finitely many vertical and horizontal subdivisions. 
See Figure \ref{Fig_patterns}. 

Next, by using patterns, we define partitions of $\Ca^2$ inductively. 
The trivial pattern corresponds to $\Ca^2$ itself. 
Assume that a rectangle in a pattern corresponds to a subset $\{w\zeta \mid \zeta \in \Ca\} \times \{w^\prime\zeta \mid \zeta \in \Ca\} \subset \Ca^2$ where $w$ and $w^\prime$ are finite words on $\{0, 1\}$. 
For the vertical subdivision, we defined that the left rectangle corresponds to the subset $\{w0\zeta \mid \zeta \in \Ca\} \times \{w^\prime\zeta \mid \zeta \in \Ca\}$ and the right rectangle corresponds to the subset $\{w1\zeta \mid \zeta \in \Ca\} \times \{w^\prime\zeta \mid \zeta \in \Ca\}$. 
Similarly, for the horizontal subdivision, we define that the bottom rectangle corresponds to the subset $\{w\zeta \mid \zeta \in \Ca\} \times \{w^\prime0\zeta \mid \zeta \in \Ca\}$ and the top rectangle corresponds to the subset $\{w\zeta \mid \zeta \in \Ca\} \times \{w^\prime1\zeta \mid \zeta \in \Ca\}$. 
Consequently, a set of rectangles of a pattern gives a partition of $\Ca^2$. 
For example, the pattern once divided vertically illustrated in Figure \ref{Fig_patterns} corresponds to $\{0\zeta\mid \zeta \in \Ca\} \times \Ca \cup \{1\zeta\mid \zeta \in \Ca\}\times \Ca$. 
We often identify a partition of $\Ca^2$ with the corresponding partition of $[0, 1]^2$. 

For a pattern with $m$ rectangles, we assign a number from $0$ to $m-1$ to each pattern. 
Such a pattern is called \textit{numbered pattern}. 
Let $P$ and $P^\prime$ be numbered patterns with $m$ rectangles. 
Let $R_i$ (resp.~ $R_i^\prime$) be a rectangle of $P$ (resp.~$P^\prime$) numbered $i$. 
Then there exist two subsets $\{a_i\zeta \mid \zeta \in \Ca\}\times \{b_i\zeta \mid \zeta \in \Ca\}$ and $\{a_i^\prime\zeta \mid \zeta \in \Ca\}\times \{b_i^\prime\zeta \mid \zeta \in \Ca\}$ corresponding to $R_i$ and $R_i^\prime$, respectively. 
By mapping each element $(a_i\zeta, b_i\zeta) \in \{a_i\zeta \mid \zeta \in \Ca\}\times \{b_i\zeta \mid \zeta \in \Ca\}$ to $(a_i^\prime\zeta, b_i^\prime\zeta) \in \{a_i^\prime\zeta \mid \zeta \in \Ca\}\times \{b_i^\prime\zeta \mid \zeta \in \Ca\}$, we obtain a homeomorphism on $\Ca^2$. 
\begin{definition}
\textit{Higher dimensional Thompson's group} $2V$ is a subgroup of $\mathop{\mathrm{Homeo}}(\Ca^2)$, where $2V$ consists of all homeomorphisms obtained from pairs of numbered patterns with the same number of rectangles. 
\end{definition}
Following a familiar conversion, we write $fg$ for $g \circ f$; namely, we always consider the right action of $2V$ on $\Ca^2$. 

%We call numbered patterns that determine the partition of the domain set of $\Ca^2$ and range set of $\Ca^2$ the \textit{domain pattern} and the \textit{range pattern} of the pair of numbered patterns, respectively. 
The \textit{domain} and \textit{range pattern} of a pair of numbered patterns are defined as the patterns that determines the partition of the domain and range set of $\Ca^2$, respectively. 
\begin{remark}
The group $nV$ is defined similarly as a subgroup of $\mathop{\mathrm{Homeo}}(\Ca^n)$. 
Using the unit $n$-cube instead of the unit square, $n$ subdivisions are defined, which yield partitions of $\Ca^n$. 
\end{remark}
Note that two distinct pairs of numbered patterns may give the same map. 
Let $(P, P^\prime)$ be a pair of numbered patterns with the same number of rectangles (such a pair is just called a \textit{pair of numbered patterns}). 
Let $R_i$ and $R_i^\prime$ be rectangles of $P$ and $P^\prime$ numbered $i$, respectively. 
Apply a vertical subdivision to $R_i$, and assign $i_1$ to the left rectangle and $i_2$ to the right rectangle. 
Do the same to $R_i^\prime$ and assign $i_1$ and $i_2$. 
Consequently, the obtained maps from the two pairs are the same. 
The same also holds for the case of horizontal subdivisions. 
A \textit{reduced} pair of numbered patterns is a pair of numbered patterns where the inverse operations (called \textit{vertical} and \textit{horizontal reductions}) can not be applied to any rectangles. 

In order to multiply two pairs of numbered patterns $(P_+, P_-)$ and $(Q_+, Q_-)$, take a common refinement $P$ of $P_-$ and $Q_+$, and by using the operations, construct pairs $(P_+^\prime, P)$ and $(P, Q_-^\prime)$ such that they give the same maps, respectively. 
Then the pair $(P_+^\prime, Q_-^\prime)$ is the desired one. 

Unfortunately, unlike Thompson's groups, there is no known way to define a unique reduced numbered pair for each element in $2V$. 
The notion of grid diagrams provides a solution to this issue. 
A detailed explanation will be given in Section \ref{Section_grid}. 
%%%%%%%%%%%%%%%%%%%%%%%%%%%%%%%%%%%%%%%%%
\subsubsection{Pairs of colored binary trees}
%%%%%%%%%%%%%%%%%%%%%%%%%%%%%%%%%%%%%%%%%
Just as elements of Thompson's groups are represented by pairs of binary trees, there exists a similar approach for $2V$. 
Like numbered patterns, no known unique representatives for elements of $2V$. 
However, because it can represent elements more simply than numbered patterns, we will frequently use colored binary trees in this paper. 

\begin{figure}[tbp]
\begin{center}
\includegraphics[width=0.4\linewidth]{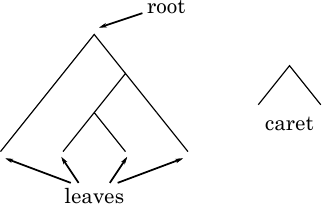}
\end{center}
\caption{A binary tree and a caret. }
\label{Fig_binary_trees}
\end{figure}

We always assume that binary trees are rooted; namely, they have one specific vertex called the \textit{root}. 
Vertices whose degree is one are called \textit{leaves}. 
A graph with only one vertex (and no edges) is also regarded as a binary tree. 
We define a \textit{caret} to be the binary tree consisting of three vertices, where the degree of the root is two and the remaining vertices are one. 
See Figure \ref{Fig_binary_trees}.

Any binary tree can be constructed inductively by attaching carets to leaves. 
A \textit{colored binary tree} is a binary tree where each caret is colored. 
We use two colors $\{a, b\}$ for $2V$ and use $n$ colors for $nV$. 
We first explain the relationship between colored binary trees and patterns. 

As previously stated, any pattern can be obtained inductively from the trivial pattern. 
Similarly, corresponding binary trees are also defined inductively. 
In the case of the trivial pattern, the corresponding binary tree is the one that consists of only the root. 
More precisely, consider the root as a leaf, the trivial pattern as a rectangle, and correspond this rectangle to this leaf. 
Since there exists no caret in this tree, we do not color it. 
Assume that for a pattern, there exists a colored binary tree with leaves such that there is a one-to-one correspondence between the set of rectangles and the set of leaves. 
Let $R$ be a rectangle of this pattern, and assume that the corresponding leaf is $i$-th (from left to right). 
For the vertical subdivision of $R$, let $R_1$ be the left rectangle and $R_2$ be the right rectangle. 
Subsequently, we attach a caret colored by $a$ to the $i$-th leaf and correspond $R_1$ to the newly created left leaf and $R_2$ to the right leaf. 
Similarly, for the horizontal subdivision, we attach a caret colored by $b$ to the $i$-th leaf and correspond the bottom rectangle to the left leaf and the top rectangle to the right leaf. 

\begin{figure}[tbp]
\begin{center}
\includegraphics[width=0.9\linewidth]{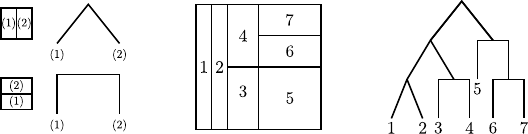}
\end{center}
\caption{The two types of carets and a colored binary tree obtained from a pattern. }
\label{Fig_colored_binary_trees}
\end{figure}

Following \cite{MR2734164}, in this paper, we represent carets colored by $a$ by the ``triangular carets'' and carets colored by $b$ by ``square carets.''
See Figure \ref{Fig_colored_binary_trees}. 
The numbers in the rectangles and under the leaves represent the one-to-one correspondence between the set of leaves and rectangles. 
Observe that numbered patterns can also be represented by colored binary trees with numbers by writing numbers under the leaves in the same way. 
Consequently, each element in $2V$ can be represented by a pair of colored binary trees with numbers. 

\begin{figure}[tbp]
\begin{center}
\includegraphics[width=0.7\linewidth]{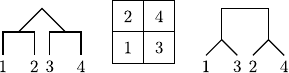}
\end{center}
\caption{Two colored binary trees give the same pattern. }
\label{Fig_quadrant}
\end{figure}

Note that, in general, more than one colored binary tree may give the same pattern. 
This is due to the fact that the patterns obtained by applying the following operations to a rectangle in a pattern are the same: the pattern obtained by subdividing vertically once and then subdividing two rectangles horizontally; subdividing horizontally once and then subdividing two rectangles vertically. 
See Figure \ref{Fig_quadrant}. 

As with pairs of numbered patterns, it is difficult to easily obtain a ``good'' pair of colored binary trees for each element of $2V$.
However, we can use this to give an estimation of the word length, as we will see in Section \ref{Section_generators}. 
For a colored binary tree $T$, we define a \textit{branch} of $T$ to be a path from the root of $T$ to a leaf of $T$. 
The \textit{depth} of $T$ is then defined as the maximum length of the branches of $T$ (with each edge having a length of one). 
It is clear that for any element $g$ in $2V$, there exists a pair of colored binary trees such that it gives $g$. 
Therefore, we may define the minimality of such pairs. 
First, we define a pair of colored binary trees to be \textit{shallow} if it is a pair of colored binary trees whose target tree has the smallest depth among pairs that give the same map. 
A pair of colored binary trees is then defined as \textit{minimal} if it is a shallow pair of binary trees with the smallest number of carets. 
It should be noted that such pairs may not be uniquely determined. 
\begin{remark}
In order to make the relationship between grid diagrams and pairs of colored binary trees clearer, the definition of minimality is slightly modified from that in \cite{MR2734164}. 
\end{remark}
%%%%%%%%%%%%%%%%%%%%%%%%%%%%%%%%%%%%%%%%%
\subsubsection{Grid diagrams} \label{Section_grid}
\begin{figure}[tbp]
\begin{center}
\includegraphics[width=0.2\linewidth]{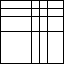}
\end{center}
\caption{A grid pattern. }
\label{Fig_grid_pattern}
\end{figure}

In this section, we explain how to represent each element in $2V$ using a ``grid'' based on \cite{burillo2024grid}.
However, unlike \cite{burillo2024grid}, the ``grid'' is constructed on the target patterns for the sake of our proof. 
Note that in \cite{burillo2024grid}, it was also pointed out that the same result holds in our setting. 
Indeed, we use only the fact that there exists a unique representative for each element in $2V$ (and $nV$). 
In this case, we only need to change the pattern we focus on from domain to target. 
Hence we omit most of the proofs in this section. 
\begin{definition}
A \textit{grid pattern} is defined as a pattern obtained by subdividing the unit square using only line segments of length one. 
%A \textit{numbered grid pattern} is a pattern with numbers where the pattern is a grid pattern. 
\end{definition}
The pattern illustrated in Figure \ref{Fig_colored_binary_trees} is not a grid pattern, since the lengths of the horizontal line segment on $[0, 1]\times\{1/2\}$ and the horizontal line segment between rectangles numbered 6 and 7 are less than one. 
An example is illustrated in Figure \ref{Fig_grid_pattern}. 
\begin{definition}
A \textit{grid diagram} is a pair of numbered patterns with the same number of rectangles, and the target pattern is a grid pattern. 
\end{definition}
It is known that any element in $2V$ can be represented by a grid diagram: 
\begin{proposition}[{\cite[Proposition 2.3]{burillo2024grid}}]
Each element in $2V$ admits a grid diagram as a representative. 
\begin{sproof}
By repeating horizontal and vertical subdivisions (such that the obtained pairs of numbered patterns always give the same map) until the target pattern, we get a grid pattern. 
\end{sproof}
\end{proposition}
Moreover, it is known that there exists a unique representative for any element in $2V$. 
In order to give such a representative, we first recall the definition that a grid diagram is reduced. 

Let $(P, G)$ be a grid diagram where $G$ is a grid pattern. 
For $G$, take a strip $I_i\times[0, 1] \subset [0, 1]\times[0, 1]$ where $I_i$ is horizontal edges of the rectangles in $G$. 
Equivalently, take a rectangle of $G$ where the rectangle is $I_i \times J_j \subset [0, 1]\times[0, 1]$, and consider $I_i \times [0, 1]$. 
We apply vertical subdivisions to all rectangles in $I_i \times [0, 1]$, and then apply vertical subdivisions to the corresponding rectangles in $P$ such that the obtained grid pattern and $(P, G)$ give the same homeomorphism on $\Ca^2$. 
We call this operation a \textit{vertical global subdivision}. 
We say a grid diagram is \textit{vertically reduced} if the inverse operation of the vertical global subdivision can not be applied to any two adjacent strips of the target pattern. 
Similarly, we can define a \textit{horizontal global subdivision} and a grid diagram to be \textit{horizontally reduced}. 
Finally, a grid diagram is said to be \textit{reduced} if it is vertically and horizontally reduced. 
Then we have the following: 
\begin{theorem}[{\cite[Theorem 3.2]{burillo2024grid}}]
Any element in $2V$ has a unique reduced grid diagram. 
\end{theorem}
Note that since each reduced grid diagram is a pair of numbered patterns, there exist corresponding pairs of colored binary trees. 
%%%%%%%%%%%%%%%%%%%%%%%%%%%%%%%%%%%%%%%%%
\subsection{A generating set of $2V$ and estimations of the word length. } \label{Section_generators}
\begin{figure}[tbp]
\begin{center}
\includegraphics[width=\linewidth]{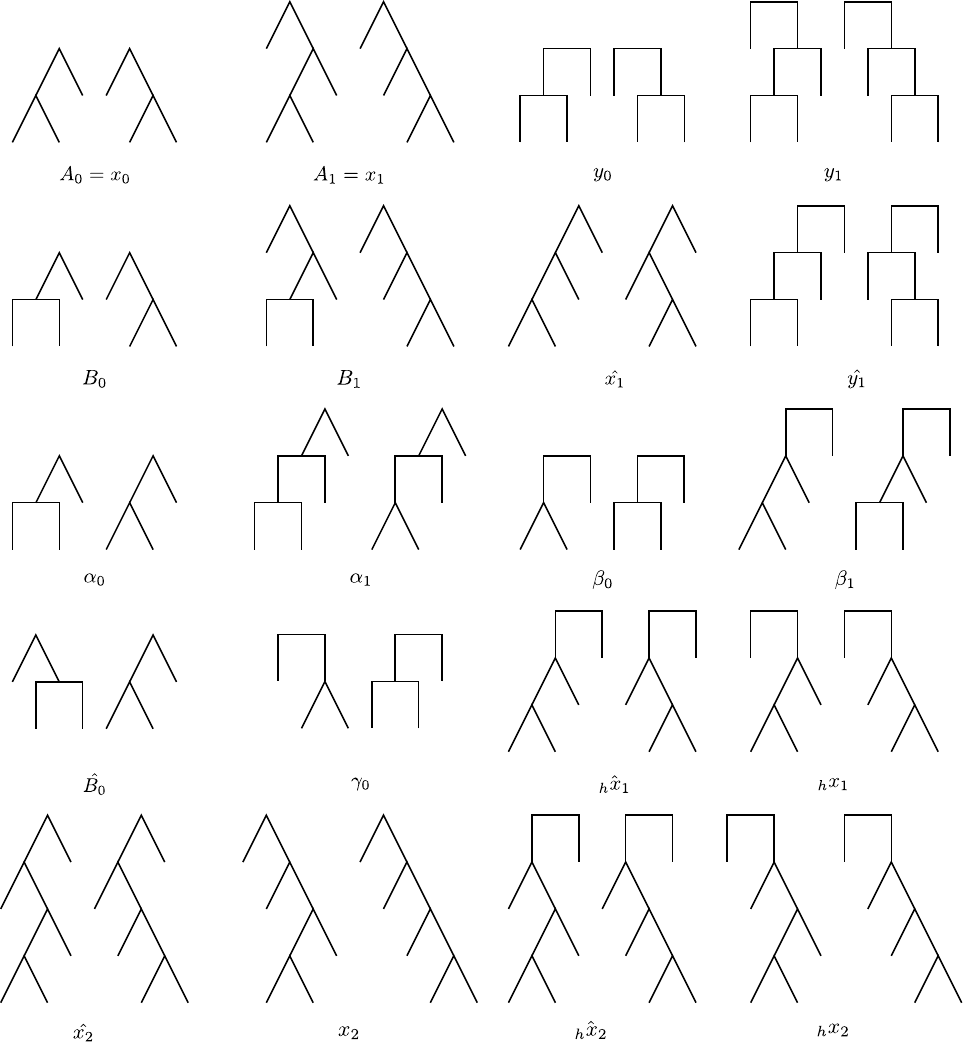}
\end{center}
\caption{Generators of $2V$ except $C_0$, $C_1$, $\pi_0$, $\pi_1$, $\overline{\pi_0}$, and $\overline{\pi_1}$. }
\label{Fig_generators}
\end{figure}

\begin{figure}[tbp]
\begin{center}
\includegraphics[width=0.8\linewidth]{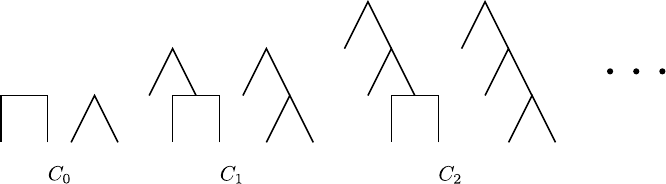}
\end{center}
\caption{Generators $C_0, C_1, \cdots$. }
\label{Fig_c0c1c2}
\end{figure}

\begin{figure}[tbp]
\begin{center}
\includegraphics[width=0.9\linewidth]{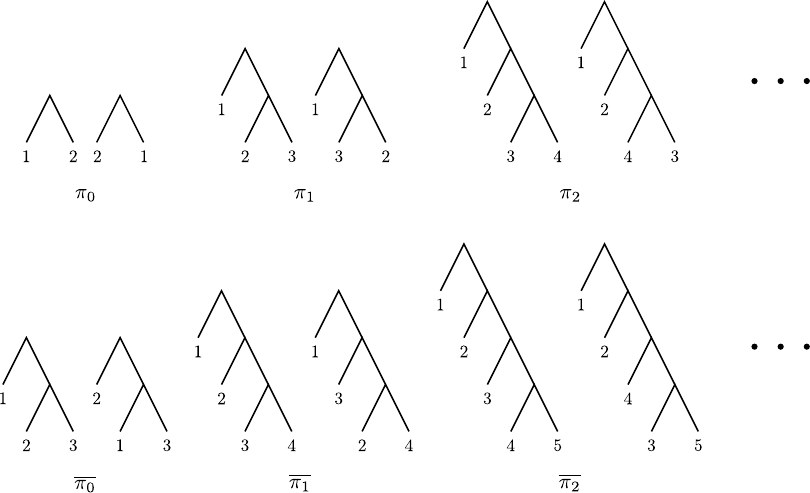}
\end{center}
\caption{Generators $\pi_0, \pi_1, \cdots$ and $\overline{\pi_0}, \overline{\pi_1}, \cdots$. }
\label{Fig_pi_i}
\end{figure}

Consider a finite set
\begin{align*}
X_{2V} \coloneqq \left\{x_0, x_1, x_2, y_i, B_i, C_i, \hat{x_j}, \hat{y_1}, \pi_i, \overline{\pi_i},  \alpha_i, \beta_i, \hat{B_0}, \gamma_0, {}_hx_j, \hat{{}_hx_j}, \;\middle|\; i \in \{0, 1\},  j\in \{1, 2\}\right\}
\end{align*}
defined in Figures \ref{Fig_generators}, \ref{Fig_c0c1c2} and \ref{Fig_pi_i}. 
The colored binary trees without numbers are all assigned $1, 2, \cdots$ to the leaves from left to right. 
Because this set contains well-known generating sets (cf.~\cite{MR2112673, MR2114568, MR2734164}), this also generates $2V$. 
Note that this set is an inefficient set specialized for calculating the divergence function of $2V$. 

Following \cite{MR2734164}, we also express an element of $2V$ as ``$P \Pi Q^{-1}$ form'' and use a part of the form in the proof of the main theorem. 
We define $A_i \coloneqq A_0^{-(i-1)}A_1A_0^{(i-1)}$ for $i\geq1$ and $B_i \coloneqq A_0^{-(i-1)}A_1A_0^{(i-1)}$ for $i\geq1$. 
Then the following holds: 
\begin{theorem}[{\cite[Theorem 2.2]{MR2734164}}]
For an element $g$ of $2V$, we have its expression $P\Pi Q^{-1}$, where
\begin{enumerate}
\item $P$ and $Q$ are represented by the following type of words: 
\begin{align*}
C_{m_1}\cdots C_{m_p} W_{i_1}\cdots W_{i_r}
\end{align*}
where $W_i$ are words on $\{A_i, B_i\}$ without inverse elements, $m_1<m_2<\cdots<m_p$, and $i_1<i_2<\cdots<i_r$. 
\item $\Pi$ are represented by a word on $\{\pi_0, \pi_1, \cdots\} \cup \{\overline{\pi_0}, \overline{\pi_1}, \cdots\}$. 
\end{enumerate}
\end{theorem}
In \cite{MR2734164}, they constructed this expression for each pair of colored binary trees. 
We also use this construction. 
See the proof of \cite[Theorem 2.1]{MR2734164} for details. 

We will consider the word metric with respect to the above generating set. 
We first recall a known result of an estimation of the word length for pairs of binary trees. 
\begin{proposition}[cf.~{\cite[Lemma 4.1]{MR2734164}}]\label{proposition_Burillo_metric}
For an element of $2V$ which is represented by a minimal pair of colored binary trees with depth $D$, its word length with respect to $X_{2V}$ is at least $D/4$. 
\end{proposition}
\begin{remark}
The difference between the denominators in \cite[Lemma 4.1]{MR2734164} and this lemma arises from the difference in the generating sets. 
The definition of minimality is also different, yet this lemma can still be shown similarly: 
if we take one of the shortest words with length $n$, we can say that the depth is at most $4n$ since each process of multiplying a generator increases the length of the branches by at most four. 
\end{remark}
As previously stated, a minimal pair of colored binary trees may not be uniquely determined for an element in $2V$. 
Hence, we use the notion of grid diagrams. 
We define the \textit{size} of a rectangle $R=\{[a_i, a_{i+1}] \times [b_j, b_{j+1}]\}$ of a pattern as $-(\log_2(a_{i+1}-a_{i})+\log_2(b_{j+1}-b_{j}))$ and write it as $\|R\|$. 
For a pattern $P$, the \textit{fineness} of $P$ is defined as the largest size of the rectangles in $P$ and written as $\|P\|$. 
For an element $g \in 2V$, the \textit{fineness} of $g$ is defined as $\|G\|$, where $ G$ is the target numbered pattern of the reduced grid diagram of $g$. 
Note that the fineness of $g$ is the same as the depth of the target tree of a pair of colored binary trees obtained from the grid diagram of $g$. 
In the proof of the main theorem, we also use the following estimations: 
%%%%%%%%%%%%%%%%%%%%%%%%%%%%%%%%%%%%%%%%%
\begin{proposition}\label{Prop_true_estimation}
Let $g$ be in $2V$ with fineness $k$. 
Then the word length of $g$ with respect to $X_{2V}$ is at least $k/8$. 
\begin{proof}
Consider a minimal pair of colored binary trees of $g$ and let $D$ be its depth. 
Take a pair of numbered patterns $(P_1, P_2)$ corresponding to this pair of colored binary trees. 
Then, we have $\|P_2\|=D$. 
We construct a grid diagram $(P, G)$ of $g$ from $(P_1, P_2)$ without vertical and horizontal global subdivisions. 
Then, since each rectangle in $G$ has horizontal and vertical lengths of at most $1/2^D$, the fineness of $G$ is at most $2D$. 
Since the fineness of $g$ is at most $\|G\|$, we have $k \leq \|G\|$. 
By Proposition \ref{proposition_Burillo_metric}, the word length of $g$ is at least $D/4$. 
Hence, the desired result is obtained. 
\end{proof}
\end{proposition}
By a similar argument, the following also holds: 
\begin{corollary}\label{Cor_true_estimation}
Let $(P_1, P_2)$ be a pair of numbered patterns corresponding to $g \in 2V$, and $R=\{w_1x_1\zeta_1 \mid \zeta_1 \in \Ca\}\times\{w_2x_2\zeta_2 \mid \zeta_2 \in \Ca\}$ be a rectangle of $P_2$ where $x_1$ and $x_2$ are in $\{0, 1\}$. 
Assume that neither of the following conditions $(1)$ nor $(2)$ is satisfied: 
\begin{enumerate}
\item
\begin{itemize}
\item $g^{-1}(R)=\{w_1^\prime x_1\zeta_1 \mid \zeta_1 \in \Ca\}\times\{w_2^\prime \zeta_2 \mid \zeta_2 \in \Ca\}$ where $w_1^\prime$ and $w_2^\prime$ are some word on $\{0, 1\}$, and
\item $g^{-1}((w_1 \hat{x_1}\zeta_1, w_2x_2\zeta_2))=(w_1^\prime\hat{x_1}\zeta_1, w_2^\prime\zeta_2)$ holds for every $(w_1 \hat{x_1}\zeta_1, w_2x_2\zeta_2) \in \{w_1 \hat{x_1}\zeta_1 \mid \zeta_1 \in \Ca\}\times \{w_2x_2\zeta_2 \mid \zeta_2 \in \Ca\}$ where $\hat{x_1}$ is in $\{0, 1\}$ with $\hat{x_1}\neq x_1$; 
\end{itemize}
\item
\begin{itemize}
\item $g^{-1}(R)=\{w_1^\prime \zeta_1 \mid \zeta_1 \in \Ca\}\times\{w_2^\prime x_2 \zeta_2 \mid \zeta_2 \in \Ca\}$ where $w_1^\prime$ and $w_2^\prime$ are some word on $\{0, 1\}$, and
\item $g^{-1}((w_1 x_1\zeta_1, w_2\hat{x_2}\zeta_2))=(w_1^\prime \zeta_1, w_2^\prime \hat{x_2}\zeta_2)$ holds for every $(w_1 x_1 \zeta_1, w_2\hat{x_2}\zeta_2) \in \{w_1 x_1\zeta_1 \mid \zeta_1 \in \Ca\}\times \{w_2\hat{x_2}\zeta_2 \mid \zeta_2 \in \Ca\}$ where $\hat{x_2}$ is in $\{0, 1\}$ with $\hat{x_2}\neq x_2$. 
\end{itemize}
\end{enumerate}
Then, the word length of $g$ is at least $\|R\|/8$. 
\begin{proof}
From the assumption, the fineness of $g$ is at least $\|R\|$. 
By Proposition \ref{Prop_true_estimation}, we have the desired result. 
\end{proof}
\end{corollary}
This assumption means that we can not apply vertical and horizontal reductions to $R$ and  ``congruent rectangles'' adjacent to $R$. 
We say a rectangle $R$ is an \textit{essential} rectangle if it satisfies the assumption of Corollary \ref{Cor_true_estimation}. 
We call $R^v$ and $R^h$ for subsets $\{w_1 \hat{x_1}\zeta_1 \mid \zeta_1 \in \Ca\}\times \{w_2x_2\zeta_2 \mid \zeta_2 \in \Ca\}$ and $\{w_1 x_1\zeta_1 \mid \zeta_1 \in \Ca\}\times \{w_2\hat{x_2}\zeta_2 \mid \zeta_2 \in \Ca\}$, respectively, which are defined from $R=\{w_1x_1\zeta_1 \mid \zeta_1 \in \Ca\}\times\{w_2x_2\zeta_2 \mid \zeta_2 \in \Ca\}$ in the assumption of Corollary \ref{Cor_true_estimation}. 
%%%%%%%%%%%%%%%%%%%%%%%%%%%%%%%%%%%%%%%%%
\subsection{Divergence functions of finitely generated groups}
For finitely generated groups, the property of having linear divergence functions is a quasi-isometric invariant. 
Since we see asymptotic properties of functions, we introduce an equivalence relation on functions from $\mathbb{R}_{>0}$ to itself as follows: 
let $f$ and $g$ be functions from $\mathbb{R}_{>0}$ to itself. 
We first define $f \preceq g$ if there exist $A, B, C, D, E \leq 0$ such that
\begin{align*}
f(x) \leq A g(Bx+C)+Dx+E
\end{align*}
holds for all $x \in \mathbb{R}_{>0}$. 
Then we define $f \approx g$ if $f \preceq g$ and $g \preceq f$ hold. 
This is an equivalence relation on the set of functions from $\mathbb{R}_{>0}$ to itself. 
Note that all linear functions and constant functions are equivalent. 

Let $G$ be a finitely generated group with a finite generating set $X$. 
Let $\Gamma$ be the Cayley graph of $G$ with respect to $X$. 
For $\delta \in (0, 1)$, we first define the $\delta$-divergence function of $\Gamma$ as follows: let $\Omega(g_1, g_2)$ be the set of all paths connecting $g_1, g_2 \in \Gamma$. 
Let $\|\omega\|$ denote the length of the path $\omega$. 
Then for $x \in \mathbb{R}_{>0}$, define the function $\phi_\delta$ by setting
\begin{align*}
\phi_\delta(x) \coloneqq \max \Big\{ \min\{\|\omega\| \mid \omega \in \Omega(g_1, g_2) \text{ with $\omega$ avoids $B(e, \delta x)$} \} \;\Big|\; |g_1|=|g_2|=x \Big\}, 
\end{align*}
where $B(e, \delta x)$ denotes the open ball of radius $\delta x$ with a center at identity of $G$, and $|g_1|, |g_2|$ denote the lengths from identity of $G$ to $g_1, g_2$ in $\Gamma$. 
If there does not exist such a path, take $\phi_\delta(x)=\infty$. 
For each $\delta \in (0, 1)$, the equivalence class of $\phi_\delta$ is invariant under quasi-isometries. 
Hence, the $\delta$-divergence function of $G$ is well-defined as an equivalence class of functions. 
\begin{definition}
We say \textit{a group $G$ has a linear divergence function} if there exists $\delta \in (0, 1)$ such that the $\delta$-divergence function of $G$ is in the equivalence class of linear maps. 
\end{definition}
From this definition, it is clear that if the $\delta$-divergence function of $G$ is in the equivalence class of linear maps, then the $\delta^\prime$-divergence function of $G$ is also in the equivalence class of linear functions for every $0<\delta^\prime \leq\delta$. 
In fact, Dru{\c{t}}u, Mozes, and Sapir \cite{MR2584607, MR3717996} showed the following theorem, which establishes a relationship between divergence functions and a topological property of asymptotic cones. 
\begin{theorem}[\cite{MR2584607, MR3717996}]
The following are equivalent: 
\begin{enumerate}
\item $G$ has a linear divergence function, 
\item For every $\delta \in (0, 1/54)$, the function $\phi_\delta$ is in the equivalence class of linear functions, and
\item None of the asymptotic cones of $G$ has a cut point. 
\end{enumerate}
\end{theorem}
%%%%%%%%%%%%%%%%%%%%%%%%%%%%%%%%%%%%%%%%%$
%%%%%%%%%%%%%%%%%%%%%%%%%%%%%%%%%%%%%%%%%
\section{Proof of Theorem \ref{theorem_divergence_nV}} \label{section_main}
\subsection{Size of the bottom left rectangle}
For a given pattern $P$, let $R_0(P)$ be defined as the rectangle of $P$ containing the point $(0, 0) \in [0, 1]\times[0, 1]$. 
We first study the change of the word length when we multiply some of the generators in $X_{2V}$. 
\begin{lemma} \label{lemma_estimations}
Let $(P_+, P_-)$ be a pair of numbered patterns representing $g \in 2V$. 
Assume that $R_0(P_-)$ is essential and a subset of $[0, 1/4]\times[0, 1]$. 
Then the following hold: 
\begin{enumerate}
\item 
There exists a pair of numbered pattern $(P_+(gx_0^{-1}), P_-(gx_0^{-1}))$ representing $gx_0^{-1}$ such that $R_0( P_-(gx_0^{-1}))$ is essential, $R_0( P_-(gx_0^{-1})) \subset [0, 1/8]\times [0, 1]$ and $\| R_0( P_-(gx_0^{-1}))\| = \|R_0(P_-)\|+1$ hold. 
\item 
There exists a pair of numbered pattern $(P_+(g\hat{B_0}), P_-(g\hat{B_0}))$ representing $g\hat{B_0}$ such that $R_0( P_-(g\hat{B_0}))$ is essential, $R_0( P_-(g\hat{B_0})) \subset [0, 1/8]\times [0, 1]$ and $\| R_0( P_-(g\hat{B_0}))\| = \|R_0(P_-)\|+1$ hold. 
\item 
If $R_0(P_-)$ is a subset of $[0, 1/4] \times [0, 1/2]$, then there exists a pair of numbered pattern $(P_+(gC_0), P_-(gC_0))$ representing $gC_0$ such that $R_0( P_-(gC_0))$ is essential, $R_0( P_-(gC_0)) \subset [0, 1/8]\times [0, 1]$ and $\| R_0( P_-(gC_0))\| = \|R_0(P_-)\|$ hold. 
\end{enumerate}
\begin{proof}
For part (1), note that we have $R_0(P_-) \cup (R_0(P_-))^v \cup (R_0(P_-))^h \subset [0, 1/2] \times [0, 1]$ if $(R_0(P_-))^h$ is defined. 
Considering the composition of pairs of patterns. 
The rectangle $R_0(P_-)$ is unchanged even if the common refinement is taken. 
Hence, $x_0^{-1}(R_0(P_-))$ is an essential rectangle of a pair of patterns representing $gx_0^{-1}$. 
The remaining claims are clear from the definition of $x_0^{-1}$. 
Note that if $(R_0(P_-))^h$ is not defined, then it is also clear since $(x_0^{-1}(R_0(P_-)))^h$ is not defined. 
By the same argument, part (2) also follows. 

For part (3), if $R_0(P_-)$ is a subset of $[0, 1/4] \times [0, 1/4]$, then it is also followed by a similar argument of the proof of part (1). 
If $R_0(P_-)$ is $[0, a] \times [0, 1/2]$ for some $a$ with $a \leq 1/4$, then $C_0(R_0(P_-))^h$ is not defined. 
Hence, we also have the desired result. 
\end{proof}
\end{lemma}
By a similar argument, we also have the following: 
\begin{lemma}
Let $(P_+, P_-)$ be a pair of numbered patterns representing $g \in 2V$. 
Assume that $R_0(P_-)$ is essential and a subset of $[0, 1]\times[0, 1/4]$. 
Then the following hold: 
\begin{enumerate}
\item 
There exists a pair of numbered pattern $(P_+(gy_0^{-1}), P_-(gy_0^{-1}))$ representing $gy_0^{-1}$ such that $R_0( P_-(gy_0^{-1}))$ is essential, $R_0( P_-(gy_0^{-1})) \subset [0, 1]\times [0, 1/8]$ and $\| R_0( P_-(gy_0^{-1}))\| = \|R_0(P_-)\|+1$ hold. 
\item 
There exists a pair of numbered pattern $(P_+(g\gamma_0), P_-(\gamma_0))$ representing $g\gamma_0$ such that $R_0( P_-(g\gamma_0))$ is essential, $R_0( P_-(g\gamma_0)) \subset [0, 1] \times [0, 1/8]$ and $\| R_0( P_-(g\gamma_0))\| = \|R_0(P_-)\|+1$ hold. 
\item 
If $R_0(P_-)$ is a subset of $[0, 1/2] \times [0, 1/4]$, then there exists a pair of numbered pattern $(P_+(gC_0^{-1}), P_-(gC_0^{-1}))$ representing $gC_0^{-1}$ such that $R_0( P_-(gC_0^{-1}))$ is essential, $R_0( P_-(gC_0^{-1})) \subset [0, 1]\times [0, 1/8]$ and $\| R_0( P_-(gC_0^{-1}))\| = \|R_0(P_-)\|$ hold. 
\end{enumerate}
\end{lemma}
%%%%%%%%%%%%%%%%%%%%%%%%%%%%%%%%%%%%%%%%%
\subsection{Construction of the path}
In the rest of this paper, we write $|\cdot|$ for the word length of $2V$ with respect to $X_{2V}$. 
For a word $w$ and $w^\prime$, we write $w\equiv w^\prime$ when $w$ and $w^\prime$ are the same as words, and $\|w\|$ denotes the length of $w$. 
If we have $w\equiv w^\prime w^{\prime \prime}$ for some words $w$, $w^\prime$ and $w^{\prime \prime}$ (the word $w^{\prime \prime}$ may be the empty word), then $w^\prime$ is said to be a \textit{prefix} of $w$ and denoted by $w^\prime \leq w$. 
Theorem \ref{theorem_divergence_nV} immediately follows from the following proposition: 
\begin{proposition} \label{proposition_path1}
There exist constants $\delta, D$ and a positive integer $Q$ such that the following holds: 
let $g \in 2V$ with $|g| \geq 4$ holds. 
Then there exists a path of length at most $D|g|$ in the Cayley graph of $2V$ which avoids a $\delta|g|$-neighborhood of the identity and which has the initial vertex $g$ and the terminal vertex $\hat{x_1}^{-Q|g|}\hat{x_2}\hat{x_1}^{Q|g|}x_1^{-Q|g|}x_2x_1^{Q|g|}$. 

In other words, there exists a word $\omega$ on the generating set such that $\|\omega\|<D|g|$; for any prefix $\omega^\prime$ of $\omega$, we have $g\omega^\prime>\delta |g|$ such that
\begin{align*}
g \omega=\hat{x_1}^{-Q|g|}\hat{x_2}\hat{x_1}^{Q|g|}x_1^{-Q|g|}x_2x_1^{Q|g|}. 
\end{align*}
%%%%%%%%%%%%%%%%%%%%%%%%%%%%%%%%%%%%%%%%%
\begin{proof}
Let $(P_+(g), P_-(g))$ be a reduced pair of numbered patterns representing $g$ such that $R_0(P_-(g))$ is essential. 
Note that $R_0(P_-(g))$ can be made essential by changing the choice of rectangles to be reduced if necessary. 
We will define six subwords, denoted by $\omega_1, \dots, \omega_6$, and then the concatenation of these subwords, $\omega_1\cdots \omega_6$, will be the desired word, $\omega$. 
In the part of subpaths \ref{subpath1}, \ref{subpath2} and \ref{subpath3}, we may assume that $R_0(P_-(g))$ is a subset of $[0, 1/2]\times[0, 1]$. 
Indeed, since $g$ is not the identity map, if not, by replacing the ``vertical argument'' with the ``horizontal one'', it is possible to join the argument from subpath \ref{subpath4} onward. 
In this case, we need other generators defined in Figure \ref{Fig_generators}, which are not used in the following arguments. 
\begin{subpath}\label{subpath1}
We define $\omega_1$ to be one of the following: 
\begin{enumerate}[label=$(\alph*)$]
\item 
If $R_0(P_-(g))$ is a subset of $[0, 1/4] \times [0, 1]$, then define $\omega_1$ to be the empty word. 
\item If $R_0(P_-(g))$ is $[0, 1/2]\times [0, 1]$, then define $\omega_1$ to be $\hat{x_1}$. 
\item If $R_0(P_-(g))$ is $[0, 1/2]\times [0, 1/2]$ and $\alpha_0(R_0(P_-(g)))$ is an essential rectangle of $g\alpha_0$, then define $\omega_1$ to be $\alpha_0$. 
\item If $R_0(P_-(g))$ is $[0, 1/2]\times [0, 1/2]$ and $\alpha_0(R_0(P_-(g)))$ is not an essential rectangle of $g\alpha_0$, then define $\omega_1$ to be $B_0 \pi_1 x_0^{-1}$ $($see Figure $\ref{Fig_path1}$$)$. 
\item If $R_0(P_-(g))$ is $[0, 1/2]\times [0, 1/2^i]$ where $i\geq2$ and $\alpha_0(R_0(P_-(g)))$ is an essential rectangle of $g\alpha_0$, then define $\omega_1$ to be $\alpha_0$. 
\item If $R_0(P_-(g))$ is $[0, 1/2]\times [0, 1/2]$ where $i\geq2$ and $\alpha_0(R_0(P_-(g)))$ is not an essential rectangle of $g\alpha_0$, then define $\omega_1$ to be $\alpha_1$. 
\end{enumerate}
\begin{figure}[tbp]
\begin{center}
\includegraphics[width=0.5\linewidth]{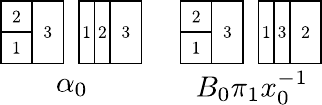}
\end{center}
\caption{Elements defined in subpath \ref{subpath1} (c) and (d). }
\label{Fig_path1}
\end{figure}
In any case, let $g_1=g \omega_1$. 
\end{subpath}
\begin{lemma} \label{lemma_subpath1}
In any case, $g_1$ has a reduced pair of numbered patterns $(P_+(g_1), P_-(g_1))$ with $R_0(P_-(g_1))$ is essential and a subset of $[0, 1/4]\times[0, 1]$. 
Moreover, $\|\omega_1\| \leq 3$ holds and for any prefix $\omega^\prime \leq \omega_1$, we have $|g\omega^\prime| \geq |g|/4$. 
\begin{proof}
%We first note that $R_0(P_-(g))$ is essential since $(P_+(g), P_-(g))$ is reduced. 
Since the latter statements are obvious, it is sufficient to show the first statement. 
But since cases (a), (c), (d), (e), and (f) of the remaining claim are also clear, we consider only case (b). 
Let $\{w_1\zeta \mid \zeta \in \Ca\} \times \{w_2 \zeta \mid \zeta \in \Ca\}$ be a rectangle corresponding to $g^{-1}(R_0(P_-(g))$. 
In the process of multiplying $\hat{x_1}$, the rectangle $g^{-1}(R_0(P_-(g))$ is subdivided into $\{w_100\zeta \mid \zeta \in \Ca\} \times \{w_2 \zeta \mid \zeta \in \Ca\}$, $\{w_101\zeta \mid \zeta \in \Ca\} \times \{w_2 \zeta \mid \zeta \in \Ca\}$, and $\{w_11\zeta \mid \zeta \in \Ca\} \times \{w_2 \zeta \mid \zeta \in \Ca\}$. 
Since the rectangles of a target pattern of $g\hat{x_1}$ corresponding to the first two are $\{00\zeta \mid \zeta \in \Ca\}\times \{\zeta \mid \zeta \in \Ca\}$ and $\{010\zeta \mid \zeta \in \Ca\}\times \{\zeta \mid \zeta \in \Ca\}$, respectively, the rectangle $\{00\zeta \mid \zeta \in \Ca\}\times \{\zeta \mid \zeta \in \Ca\}$ of $g\hat{x_1}$ is essential. 
\end{proof}
\end{lemma}
The idea of the following subpath comes from \cite{lucy2024divergence}. 
\begin{subpath}\label{subpath2}
We fix an integer $M \geq 100$. 
Consider the expression $P \Pi Q^{-1}$ of a minimal pair of colored binary trees of $g_1$ and let $C_{m_1}\cdots C_{m_p}$ be the maximal words on $\{C_i\}$ contained in $Q$, where $0 \leq m_1<\cdots <m_p$. 
When $m_1$ is not zero, we have
\begin{align*}
C_{m_1} \cdots C_{m_p}&=
(x_0^{-(m_1-1)}C_1x_0^{(m_1-1)} \cdots x_0^{-(m_p-1)}C_1x_0^{(m_p-1)}) \\ 
&=x_0^{-(m_1-1)}C_1 x_0^{-(m_2-m_1)} \cdots C_1x_0^{-(m_p-m_{p-1})}C_1x_0^{(m_p-1)} \\
&=x_0^{-(m_1-1)}\hat{B_0} x_0^{-(m_2-m_1-1)} \cdots \hat{B_0} x_0^{-(m_p-m_{p-1}-1)}\hat{B_0}x_0^{m_p}, 
%&\quad\cdot (x_0^{-M|g|}x_1x_0^{M|g|}) \cdot (x_0^{-(m_p-1)}C_1^{-1}x_0^{(m_p-1)} \cdots x_0^{-(m_1-1)}C_1^{-1}x_0^{(m_1-1)})
\end{align*}
and define
\begin{align*}
\omega_2 &\equiv x_0^{-(m_1-1)}\hat{B_0} x_0^{-(m_2-m_1-1)} \cdots \hat{B_0} x_0^{-(m_p-m_{p-1}-1)}\hat{B_0} x_0^{-(M|g_1|-m_p)}  \\
&\quad \cdot x_1 \cdot x_0^{(M|g_1|-m_p)}\hat{B_0}^{-1}x_0^{(m_p-m_{p-1}-1)} \hat{B_0}^{-1} \cdots x_0^{(m_2-m_1-1)} \hat{B_0}^{-1}x_0^{(m_1-1)}. 
\end{align*}
When $m_1$ is zero, define
\begin{align*}
\omega_2 &\equiv C_0 x_0^{-(m_2-1)}\hat{B_0} x_0^{-(m_3-m_2-1)} \cdots \hat{B_0} x_0^{-(m_p-m_{p-1}-1)}\hat{B_0} x_0^{-(M|g_1|-m_p)}  \\
&\quad \cdot x_1 \cdot x_0^{(M|g_1|-m_p)}\hat{B_0}^{-1}x_0^{(m_p-m_{p-1}-1)} \hat{B_0}^{-1} \cdots x_0^{(m_3-m_2-1)} \hat{B_0}^{-1}x_0^{(m_2-1)} C_0^{-1}. 
\end{align*}
Let $g_2=g_1 \omega_2$. 
\end{subpath}
\begin{lemma} \label{lemma_subpath2}
\begin{enumerate}
\item We have $\|\omega_2\| < 4M|g|$. 
\item For every prefix $\omega^\prime$ of $\omega_2$, we have $|g_1\omega^\prime|>|g|/64$. 
\end{enumerate}
\begin{proof}
The proof is only provided when $m_1$ is not zero; however, it can be shown similarly when $m_1$ is zero. 
We first note that $m_p \leq 4|g_1|$ holds by Proposition \ref{proposition_Burillo_metric}. 
A straightforward calculation yields the upper bound of $\|\omega_2\|$. 
Indeed, we have
\begin{align*}
&\big( (m_1-1)+1+(m_2-m_1)+1+\cdots+(m_p-m_{p-1})+1+(M|g_1|-m_p) \big) \times 2+1\\
&=2M|g_1|+1 \\
%&\leq 9|g|+36+2M|g|+1 \\
&\leq 2M(|g|+3)+1 \\
&< 4M|g|. 
\end{align*}
%Note that we have $p-1, m_p-1 \leq 3|g_1|\leq 3|g|+18$ by Proposition \ref{proposition_Burillo_metric}. 
%We can estimate similarly when $m_1=0$. 

For part (2), we first consider $\omega^\prime \leq \omega_2$ which does not contain $x_1$. 
%we first consider that $\omega^\prime$ is a prefix of $\omega_1$. 
%Since $|g|\geq \textcolor{red}{8}$ and $\|\omega_1\|\leq 6$, we have $|g\omega^\prime|>|g|-6\geq |g|/20$. 
%Next, let $\omega^\prime\equiv \omega_1\omega^{\prime \prime}$ and assume that $\omega^{\prime \prime}$ does not contain $x_1$. 
When $\|\omega^{\prime}\| \leq \lfloor |g_1|/2 \rfloor$ holds, we have $|g_1\omega^{\prime}|\geq|g_1|/2$. 
Indeed, if not, we have 
\begin{align*}
|g_1|=|g_1\omega^{\prime}(\omega^{\prime})^{-1}|\leq |g_1\omega^{\prime}|+\|\omega^{\prime}\|<\frac{|g_1|}{2}+\frac{|g_1|}{2}. 
\end{align*}
Hence we have $|g_1 \omega^\prime| \geq (|g|-3)/2 \geq |g|/8$.
Next,  we assume that $\omega^{\prime}$ does not contain $x_1$ and $\|\omega^{\prime}\| > |g_1|/2$ holds. 
Then by Lemmas \ref{lemma_subpath1} and \ref{lemma_estimations}, there exists a pair of numbered patterns $(P_+(g_1\omega^\prime), P_-(g_1\omega^\prime))$ representing $g_1\omega^\prime$ such that $R_0(P_-(g_1\omega^\prime))$ is essential. 
Then we have
\begin{align}
\|R_0(P_-(g_1\omega^\prime))\| = \|R_0(P_-(g_1))\|+\|\omega^\prime\| > \|\omega^\prime\|>\frac{|g_1|}{2}. \label{ineq_subpath2}
\end{align}
By Corollary \ref{Cor_true_estimation}, we have
\begin{align*}
|g_1\omega^\prime| \geq \frac{\|R_0(P_-(g_1\omega^\prime))\|}{8} >\frac{|g_1|}{16}\geq \frac{|g|}{64}. 
\end{align*}

Next, we assume that $\omega^\prime$ contains $x_1$ and no $\hat{B_0}^{-1}$. 
Then there exists $i \geq 0$ such that we have 
\begin{align*}
&g_1 \omega^\prime \\
&=g_1x_0^{-(m_1-1)}\hat{B_0} x_0^{-(m_2-m_1-1)} \cdots \hat{B_0} x_0^{-(m_p-m_{p-1}-1)}\hat{B_0} (x_0^{-(M|g_1|-m_p)} x_1 x_0^{(M|g_1|-m_p)}) x_0^{-i}
\end{align*}
as an element in $2V$. 
Since $x_0^{-(M|g_1|-m_p)} x_1 x_0^{(M|g_1|-m_p)}$ is identity on $[0, 1/2] \times [0, 1]$, the essentiality of the rectangle is preserved when $i=0$. 
Hence, by inequality \eqref{ineq_subpath2} and Lemma \ref{lemma_estimations}, we also have $|g_1\omega^\prime| \geq |g|/96$. 

Finally, consider the remaining case. 
From the construction of $\omega_2$, it can be observed that $g_2$ is obtained by attaching carets to the minimal pair of colored binary trees $(T_+(g_1), T_-(g_1))$. 
This implies that there exists a pair of numbered patterns of $g_2$ such that the rectangle of the target pattern that contains a $(1, 1) \in [0, 1]\times [0, 1]$ is essential, and its size is at least $M|g_1|+2-4|g_1|$ by Proposition \ref{proposition_Burillo_metric}. 
See Figure \ref{Fig_path2} for an illustration of this argument. 
\begin{figure}[tbp]
\begin{center}
\includegraphics[width=\linewidth]{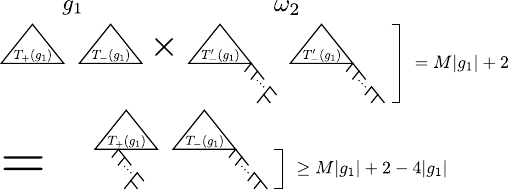}
\end{center}
\caption{An illustration of the multiplication of $g_1$ and $\omega_2$, where $(T_+(g_1), T_-(g_1))$ is a minimal pair of colored binary trees of $g_1$ and $T_-^\prime(g_1)$ is the maximal all-right colored binary tree contained $T_-(g_1)$. }
\label{Fig_path2}
\end{figure}
We now consider $g_1\omega^\prime$ as $g_2 \omega^{\prime \prime}$ by a certain word $\omega^{\prime \prime}$. 
Note that $\|\omega^{\prime \prime}\| \leq m_p$ holds. 
Since $m_p \leq 4|g_1|$ holds, by Corollary \ref{lemma_subpath1} and Lemma \ref{lemma_subpath1}, we have
\begin{align}
|g_1\omega^\prime|=
|g_2\omega^{\prime \prime}|\geq |g_2|-\|\omega^{\prime \prime}\| \geq \frac{M|g_1|+2-4|g_1|}{8}-4|g_1|>\frac{1}{8}(M-36)|g_1|\geq 8|g_1| \geq 2|g|. \label{g2_estimation}
\end{align}
Hence for any prefix $\omega^\prime \leq \omega_2$, we have $|g_1\omega^\prime|>|g|/64$. 
\end{proof}
\end{lemma}
\begin{subpath} \label{subpath3}
Let $\omega_3$ be a minimal word on $X_{2V}$ such that $\omega_3=g_1^{-1}$ holds. 
Let $g_3=g_2\omega_3$. 
\end{subpath}
\begin{lemma} \label{lemma_subpath3}
\begin{enumerate}
\item We have $\|\omega_3\| \leq 2|g|$. 
\item For every prefix $\omega^\prime$ of $\omega_3$, we have $|g_2 \omega^\prime| > 2|g|$. 
\end{enumerate}
\begin{proof}
The first claim is obvious by Lemma \ref{lemma_subpath1}. 
For the second claim, we note that $|\omega^\prime| \leq |g_1|$ holds since $\omega_3$ is a minimal word. 
Then we have
\begin{align*}
|g_2 \omega^\prime|\geq  |g_2|- |\omega^\prime| \geq |g_2|-|g_1|>2|g|, 
\end{align*}
as estimated in inequality \eqref{g2_estimation}. 
\end{proof}
\end{lemma}
\begin{subpath} \label{subpath4}
We fix an integer $Q \geq 48M$. 
We define $\omega_4$ based on a subset $[0, 1] \times [0, 1]$ where $g_3$ is the identity map. 
From the construction, one of the following holds (see also Figure $\ref{Fig_g3}$): 
\begin{figure}[tbp]
\begin{center}
\includegraphics[width=0.7\linewidth]{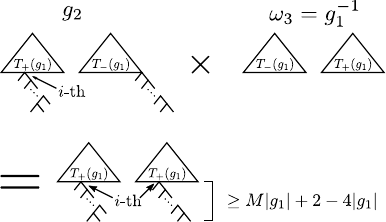}
\end{center}
\caption{An illustration of $g_3$. Observe that $g_3$ is the identity map on the subset corresponding to all leaves except the $i$-th  leaf. }
\label{Fig_g3}
\end{figure}
\begin{enumerate}[label=$(\alph*)$]
\item $g_3$ is the identity map on $[0, 1] \times [0, 1/2]$; 
\item $g_3$ is the identity map on $[0, 1] \times [1/2, 1]$; 
\item $g_3$ is the identity map on $[0, 1/2] \times [0, 1]$; and
\item $g_3$ is the identity map on $[1/2, 1] \times [0, 1]$. 
\end{enumerate}
Then we define $\omega_4$ as the word determined from one of the following cases of capital letters corresponding to each small letter: 
\begin{enumerate}[label=$(\Alph*)$]
\item let $\omega_4 \equiv \hat{{}_hx_1}^{-(Q|g|+1)}\hat{{}_hx_2}\hat{{}_hx_1}^{Q|g|+1}$; 
\item let $\omega_4 \equiv {{}_hx_1}^{-(Q|g|+1)} {}_hx_2 {{}_hx_1}^{Q|g|+1}$; 
\item let $\omega_4 \equiv \hat{x_1}^{-Q|g|}\hat{x_2}\hat{x_1}^{Q|g|}$; and 
\item let $\omega_4 \equiv {x_1}^{-Q|g|} {x_2} {x_1}^{Q|g|}$. 
\end{enumerate}
Let $g_4=g_3 \omega_4$. 
\end{subpath}
\begin{lemma} \label{lemma_subpath4}
\begin{enumerate}
\item We have $\| \omega_4\| \leq 3Q|g|$. 
\item For every prefix $\omega^\prime$ of $\omega_4$, we have $|g_3\omega^\prime| > 3|g|$. 
\item As elements in $2V$, $g_3$ and $\omega_4$ commute. 
\end{enumerate}
\begin{proof}
Part (1) follows from a straightforward estimation, and part (3) is obvious since the supports of $g_3$ and $\omega_4$ are disjoint. 
For part (2), in any case, we note that generators in $\omega_4$ preserve the rectangle with its size at least $M|g_1|+2-4|g_1|$. 
Since the process of obtaining the essential rectangle from this rectangle requires only at most one horizontal reduction, $g_3\omega_4$ is also represented by a pair of numbered patterns with an essential rectangle of size at least $M|g_1|+2-4|g_1|$. 
Hence, by Corollary \ref{Cor_true_estimation} and Lemma \ref{lemma_subpath1}, we have
\begin{align*}
|g_3\omega^\prime| \geq \frac{M|g_1|+2-4|g_1|}{8}>12|g_1| \geq 3|g|, 
\end{align*}
which is the desired result. 
\end{proof}
\end{lemma}
\begin{subpath} \label{subpath5}
Let $\omega_5$ be a minimal word on $X_{2V}$ such that $\omega_5=g_3^{-1}$ holds. Let $g_5=g_4\omega_5$. 
\end{subpath}
\begin{lemma} \label{lemma_subpath5}
\begin{enumerate}
\item We have $\| \omega_5\| \leq 5M|g|$. 
\item For every prefix $\omega^\prime$ of $\omega_5$, we have $|g_4 \omega^\prime| > M|g|$. 
\end{enumerate}
\begin{proof}
For part (1), by Lemmas \ref{lemma_subpath1}, \ref{lemma_subpath2} and \ref{lemma_subpath3}, we have
\begin{align*}
\|\omega_5\| \leq |g|+ \|\omega_1\|+\|\omega_2\|+\|\omega_3\| \leq |g|+3+4M|g|+2|g| <5M|g|. 
\end{align*}

For part (2), we first note that $\omega_4$ is represented by a pair of patterns with an essential rectangle of size at least $Q|g|+3$. 
In particular, by a similar argument to the proof of Lemma \ref{lemma_subpath4}, the horizontal length of this rectangle is unchanged for $g_4=g_3\omega_4$. 
Hence by Corollary \ref{Cor_true_estimation}, we have $|g_4|\geq (Q|g|+2)/8$. 
Therefore we have
\begin{align*}
|g_4 \omega^\prime| \geq |g_4|-\|\omega_5\| > 6M|g|-5M|g|=M|g|. 
\end{align*}
This completes the proof. 
\end{proof}
\end{lemma}
We may now obtain $g\omega_1 \omega_2 \omega_3 \omega_4 \omega_5=\omega_4$ which only depends on $|g|$. 
The final step is to connect any of cases (A) to (D) defined in subpath \ref{subpath4} to $\hat{x_1}^{-Q|g|}\hat{x_2}\hat{x_1}^{Q|g|}x_1^{-Q|g|}x_2x_1^{Q|g|}$ by a final subpath. 
In order to define this subpath, we write the subpaths defined in cases (A) to (D) as $\omega_4(A)$, $\omega_4(B)$, $\omega_4(C)$ and $\omega_4(D)$, respectively. 
\begin{subpath} \label{subpath6}
If the path $\omega_4$ is $\omega_4(A)$, let $\omega_6 \equiv \omega_4(B) \omega_4(C)$. 
If the path $\omega_4$ is $\omega_4(B)$, let $\omega_6 \equiv \omega_4(A) \omega_4(C)$. 
If the path $\omega_4$ is $\omega_4(C)$, let $\omega_6 \equiv \omega_4(D)$. 
Finally, if the path $\omega_4$ is $\omega_4(D)$, let $\omega_6 \equiv \omega_4(C)$. 
\end{subpath}
\begin{lemma} \label{lemma_subpath6}
\begin{enumerate}
\item In any case, we have $g_5\omega_6=\hat{x_1}^{-Q|g|}\hat{x_2}\hat{x_1}^{Q|g|}x_1^{-Q|g|}x_2x_1^{Q|g|}$ as an element in $2V$. 
\item We have $\|\omega_6\| \leq 6Q|g|$. 
\item For every prefix $\omega^\prime$ of $\omega_6$, we have $|g_5 \omega^\prime| > 6M|g|$. 
\end{enumerate}
\begin{proof}
For part (1), we note that $\omega_4(A)\omega_4(B)=\omega_4(D)$ holds as elements in $2V$. 
Observe that the supports of $\omega_4(A)$ and $\omega_4(B)$ are disjoint, and the same holds for $\omega_4(C)$ and $\omega_4(D)$. 
Since we have $\omega_4(C)\omega_4(D)=\hat{x_1}^{-Q|g|}\hat{x_2}\hat{x_1}^{Q|g|}x_1^{-Q|g|}x_2x_1^{Q|g|}$, we obtain the desired result. 
Part (2) follows from Lemma \ref{lemma_subpath4}. 

Finally, part (3) follows from the following observation: in the process of multiplying generators of $\omega_6$, there always exists a pair of numbered patterns with an essential rectangle whose horizontal length is at least $Q|g|+2$. 
Hence by Corollary \ref{Cor_true_estimation}, we have $|g_5\omega^\prime|\geq (Q|g|+2)/8>6M|g|$. 
\end{proof}
\end{lemma}
Now, we define $\omega$ as $\omega_1\cdots \omega_6$. 
From Lemmas \ref{lemma_subpath1}, \ref{lemma_subpath2}, \ref{lemma_subpath3}, \ref{lemma_subpath4}, \ref{lemma_subpath5}, \ref{lemma_subpath6}, take $D=10Q$ and $\delta=1/64$. 
Then we have that $g\omega=\hat{x_1}^{-Q|g|}\hat{x_2}\hat{x_1}^{Q|g|}x_1^{-Q|g|}x_2x_1^{Q|g|}$, $\|\omega\| < D|g|$, and $|g\omega^\prime|>\delta|g|$ for any prefix $\omega^\prime$ of $\omega$. 
This completes the proof of Proposition \ref{proposition_path1}. 
%%%%%%%%%%%%%%%%%%%%%%%%%%%%%%%%%%%%%%%%%$
\end{proof} \end{proposition} %end of the proof of Proposition \ref{proposition_path1}
%%%%%%%%%%%%%%%%%%%%%%%%%%%%%%%%%%%%%%%%%$
%%%%%%%%%%%%%%%%%%%%%%%%%%%%%%%%%%%%%%%%%

\section*{Acknowledgements}
The author would like to thank Professor Tomohiro Fukaya for his helpful comments. 
The author also appreciates Professor Lucy Lifschitz for teaching about the divergence function of the Lodha--Moore group. 
%%%%%%%%%%%%%%%%%%%%%%%%%%%%%%%%%%%%%%%%%
%%%%%%%%%%%%%%%%%%%%%%%%%%%%%%%%%%%%%%%%%
\bibliographystyle{plain}
\bibliography{references} 

\begin{bibdiv}
\begin{biblist}

\bib{bleak2010family}{article}{
      author={Bleak, Collin},
      author={Lanoue, Daniel},
       title={A family of non-isomorphism results},
        date={2010},
     journal={Geom. Dedicata},
      volume={146},
      number={1},
       pages={21\ndash 26},
}

\bib{bridson2013metric}{book}{
      author={Bridson, Martin~R.},
      author={Haefliger, Andr{\'e}},
       title={Metric spaces of non-positive curvature},
      series={Grundlehren Math. Wiss.},
   publisher={Berlin: Springer},
        date={1999},
      volume={319},
        ISBN={3-540-64324-9},
}

\bib{MR2112673}{article}{
      author={Brin, Matthew~G.},
       title={Higher dimensional {T}hompson groups},
        date={2004},
        ISSN={0046-5755,1572-9168},
     journal={Geom. Dedicata},
      volume={108},
       pages={163\ndash 192},
         url={https://doi.org/10.1007/s10711-004-8122-9},
}

\bib{MR2114568}{article}{
      author={Brin, Matthew~G.},
       title={Presentations of higher dimensional {T}hompson groups},
        date={2005},
        ISSN={0021-8693,1090-266X},
     journal={J. Algebra},
      volume={284},
      number={2},
       pages={520\ndash 558},
         url={https://doi.org/10.1016/j.jalgebra.2004.10.028},
}

\bib{brin2010baker}{article}{
      author={Brin, Matthew~G.},
       title={{On the baker's map and the simplicity of the higher dimensional
  Thompson groups $nV$}},
        date={2010},
     journal={Publ. Math.},
       pages={433\ndash 439},
}

\bib{MR2734164}{article}{
      author={Burillo, Jos\'{e}},
      author={Cleary, Sean},
       title={Metric properties of higher-dimensional {T}hompson's groups},
        date={2010},
        ISSN={0030-8730,1945-5844},
     journal={Pacific J. Math.},
      volume={248},
      number={1},
       pages={49\ndash 62},
         url={https://doi.org/10.2140/pjm.2010.248.49},
}

\bib{burillo2024grid}{article}{
      author={Burillo, José},
      author={Cleary, Sean},
      author={Nucinkis, Brita},
       title={{Grid diagrams for higher-dimensional Thompson's groups}},
        date={2024},
     journal={arXiv preprint arXiv:2403.02562},
}

\bib{cannon1996introductory}{article}{
      author={Cannon, J.~W.},
      author={Floyd, W.~J.},
      author={Parry, W.~R.},
       title={Introductory notes on {R}ichard {T}hompson's groups},
        date={1996},
        ISSN={0013-8584},
     journal={Enseign. Math. (2)},
      volume={42},
      number={3-4},
       pages={215\ndash 256},
}

\bib{MR2584607}{article}{
      author={Dru\c{t}u, Cornelia},
      author={Mozes, Shahar},
      author={Sapir, Mark},
       title={Divergence in lattices in semisimple {L}ie groups and graphs of
  groups},
        date={2010},
        ISSN={0002-9947,1088-6850},
     journal={Trans. Amer. Math. Soc.},
      volume={362},
      number={5},
       pages={2451\ndash 2505},
         url={https://doi.org/10.1090/S0002-9947-09-04882-X},
}

\bib{MR3717996}{article}{
      author={Dru\c{t}u, Cornelia},
      author={Mozes, Shahar},
      author={Sapir, Mark},
       title={Corrigendum to ``{D}ivergence in lattices in semisimple {L}ie
  groups and graphs of groups''},
        date={2018},
        ISSN={0002-9947,1088-6850},
     journal={Trans. Amer. Math. Soc.},
      volume={370},
      number={1},
       pages={749\ndash 754},
         url={https://doi.org/10.1090/tran/7376},
}

\bib{gersten1994quadratic}{article}{
      author={Gersten, Stephen~M},
       title={Quadratic divergence of geodesics in {CAT(0)} spaces},
        date={1994},
     journal={Geom. Funct. Anal.},
      volume={4},
      number={1},
       pages={37\ndash 51},
}

\bib{golan2019divergence}{article}{
      author={Golan, Gili},
      author={Sapir, Mark},
       title={Divergence functions of {T}hompson groups},
        date={2019},
        ISSN={0046-5755},
     journal={Geom. Dedicata},
      volume={201},
       pages={227\ndash 242},
         url={https://doi-org.utokyo.idm.oclc.org/10.1007/s10711-018-0390-x},
}

\bib{zbMATH00437296}{book}{
      author={Gromov, Mikhael},
       title={Geometric group theory, {Volume} 2: {Asymptotic} invariants of
  infinite groups},
      series={Lond. Math. Soc. Lect. Note Ser.},
   publisher={Cambridge: Cambridge University Press},
        date={1993},
      volume={182},
        ISBN={0-521-44680-5},
}

\bib{hennig2012presentations}{article}{
      author={Hennig, Johanna},
      author={Matucci, Francesco},
       title={{Presentations for the higher-dimensional Thompson groups $nV$}},
        date={2012},
     journal={Pacific J. Math.},
      volume={257},
      number={1},
       pages={53\ndash 74},
}

\bib{issini2023linear}{article}{
      author={Issini, Letizia},
       title={On linear divergence in finitely generated groups},
        date={2023},
     journal={arXiv preprint arXiv:2311.12938},
}

\bib{kodama2023divergence}{article}{
      author={Kodama, Yuya},
       title={Divergence function of the braided {Thompson} group},
        date={2023},
        ISSN={2156-2261},
     journal={Kyoto J. Math.},
      volume={63},
      number={2},
       pages={435\ndash 470},
  url={projecteuclid.org/journals/kyoto-journal-of-mathematics/volume-63/issue-2/Divergence-function-of-the-braided-Thompson-group/10.1215/21562261-10428532.full},
}

\bib{lucy2024divergence}{article}{
      author={Lifschitz, Lucy},
       title={{Divergence of the Lodha--Moore group}},
        date={in press},
     journal={Proc. Am. Math. Soc.},
      eprint={https://doi.org/10.1090/proc/14987},
}

\bib{sheng2024divergence}{article}{
      author={Sheng, Xiaobing},
       title={{Divergence} {Property} of the {Brown--Thompson} {Groups} and
  {Braided} {Thompson} {Groups}},
        date={online first},
     journal={Transform. Groups},
      eprint={https://doi.org/10.1007/s00031-023-09839-8},
}

\end{biblist}
\end{bibdiv}
%%%%%%%%%%%%%%%%%%%%%%%%%%%%%%%%%%%%%%%%%
%%%%%%%%%%%%%%%%%%%%%%%%%%%%%%%%%%%%%%%%%
%%%%%%%%%%%%%%%%%%%%%%%%%%%%%%%%%%%%%%%%%

\bigskip
Yuya Kodama

\address{
Graduate School of Science and Engineering, Kagoshima University, 
1-21-35 Korimoto, Kagoshima-city, Kagoshima, 890-0065, Japan
}

\textit{E-mail address}: \href{mailto:yuya@sci.kagoshima-u.ac.jp}{\texttt{yuya@sci.kagoshima-u.ac.jp}}
%%%%%%%%%%%%%%%%%%%%%%%%%%%%%%%%%%%%
\end{document}